\documentclass[12pt]{amsart}
\usepackage{amssymb}

\numberwithin{equation}{section}

\newcommand{\fullref}[1]{\ref{#1} on page~\pageref{#1}}

\newcommand{\ndash}{\nobreakdash-\hspace{0pt}}
\newcommand{\Ndash}{\nobreakdash--}

\newcommand{\dd}{{\mathrm{d}}}

\DeclareMathOperator{\rank}{rank}

\newcommand{\toto}{\rightrightarrows}

\newcommand{\frX}{\mathfrak{X}}

\newcommand{\pish}{\pi^\#}

\newcommand{\gtwo}{G^{(2)}}
\newcommand{\wA}{{\widetilde A}}
\newcommand{\wB}{{\widetilde B}}
\newcommand{\wG}{{\widetilde G}}
\newcommand{\wOmega}{{\widetilde \Omega}}
\newcommand{\womega}{{\widetilde \omega}}
\newcommand{\hOmega}{{\widehat \Omega}}
\newcommand{\homega}{{\widehat \omega}}

\newtheorem{Thm}{Theorem}[section]
\newtheorem{Prop}[Thm]{Proposition}
\newtheorem{Lem}[Thm]{Lemma}

\newtheorem*{Thm*}{Theorem}

\theoremstyle{remark}

\newtheorem*{Ack}{Acknowledgment}

\theoremstyle{definition}

\newcommand{\braket}[2]{\left\langle{\,{#1}\,,\,{#2}\,}\right\rangle}

\newcommand{\Lie}[2]{{\left[{\,{#1}\,,\,{#2}\,}\right]}}

\newcommand{\Poiss}[2]{\left\{{\,{#1}\,,\,{#2}\,}\right\}}

\newcommand{\LIE}{\mathsf{Lie}}
\newcommand{\G}{\mathsf{G}}
\newcommand{\PP}{\mathsf{P}}

\newcommand{\de}{\partial}

\newcommand{\calC}{\mathcal{C}}

\newcommand{\calL}{\mathcal{L}}

\newcommand{\calP}{\mathcal{P}}

\newcommand{\LL}{{\mathrm{L}}}

\def\gpd{\,\lower1pt\hbox{$\longrightarrow$}\hskip-.24in\raise2pt
               \hbox{$\longrightarrow$}\,}

\newcommand{\be}{\begin{eqnarray*}}
\newcommand{\ee}{\end{eqnarray*}}

\newcommand\qq{\em}
\newcommand\cmp[1]{{\qq Commun.\ Math.\ Phys.\ \bf #1}}

\newcommand\anm[1]{{\qq Ann.\ Math.\ \bf #1}}

\newcommand\ptps[1]{{\qq Prog.\ Theor.\ Phys.\ Suppl.\ \bf #1}}
\newcommand\ijm[1]{{\qq Int.\ J. Math.\ \bf #1}}
\newcommand\dmj[1]{{\qq Duke Math.\ J. \bf #1}}
\newcommand\jgp[1]{{\qq J. Geom.\ Phys.\ \bf #1}}
\newcommand\ajm[1]{{\qq Am.\ J. Math.\ \bf #1}}
\newcommand\jmsj[1]{{\qq J. Math.\ Soc.\ Japan \bf #1}}

\begin{document} 

\title[Integration of Poisson manifolds and Lie algebroids]
{On the integration of Poisson manifolds, Lie algebroids, and coisotropic submanifolds}
\date{August 27, 2003}

\author[A.~S.~Cattaneo]{Alberto~S.~Cattaneo}
\address{Institut f\"ur Mathematik, Universit\"at Z\"urich--Irchel,  
Winterthurerstrasse 190, CH-8057 Z\"urich, Switzerland}  
\email{asc@math.unizh.ch}

\thanks{A.~S.~C. acknowledges partial support of SNF Grant No.~20-100029/1}

\begin{abstract}
In recent years methods for the integration of Poisson manifolds and of Lie algebroids
have been proposed, the latter being usually presented as a generalization of the former.
In this note it is shown that the latter method is actually related to (and may be derived from)
a particular case of
the former if one regards dual of Lie algebroids as special Poisson manifolds. The core
of the proof is the fact, discussed in the second part of this note, that coisotropic submanifolds
of a (twisted) Poisson manifold are in one-to-one correspondence with possibly singular 
Lagrangian subgroupoids
of source-simply-connected (twisted) symplectic groupoids.
\end{abstract}

\maketitle

\section{Introduction}\label{sec-intro}
The ``infinitesimal form'' $\LIE(G)$
of a Lie groupoid $G$ (viz., the bundle of vectors tangent to the source-fibers 
restricted to the manifold of units) is naturally equipped with a Lie algebroid structure
(we will recall definitions and basic facts in Sect.~\ref{sec-prel}).
On the other hand, there exist Lie algebroids that do not arise this way. The special ones in the image
of the $\LIE$\ndash functor are called integrable. When a Lie algebroid is integrable, then there is a unique
source-simply-connected ({\sf ssc}) Lie groupoid, up to isomorphisms \cite{MM}.

The cotangent bundle of a Poisson manifold may be given the structure of a Lie algebroid with anchor map
induced from the Poisson bivector field and the Lie bracket on $1$\ndash forms given by 
Koszul \cite{Kosz}.
When this Lie algebroid is integrable, one says that the Poisson manifold is integrable.
The ssc Lie groupoid can in this case be endowed
with a multiplicative symplectic form and is called a symplectic groupoid. 
The basic example is the cotangent
bundle $T^*G$ of a Lie group $G$ as a Lie groupoid over the Poisson manifold $(\LIE(G))^*$.

As we will recall in Sect.~\ref{sec-ipm},
a method was introduced in \cite{CaFeOb} to integrate Poisson manifolds to symplectic groupoids
by symplectic reduction
from an infinite-\hspace{0pt}dimensional manifold (the cotangent bundle of the path space
of the Poisson manifold). 
The symplectic quotient is in general singular, and actually this
happens if{f} the Poisson manifold is not integrable.
The authors of
\cite{S} and \cite{CF} independently observed that the above method allows for a generalization to any 
Lie algebroid.
The main contribution of \cite{CF} is then to use this construction to characterize integrable Lie algebroids
in terms of an if-and-only-if criterion.
See also \cite{CFPoiss} for a general discussion of obstructions to 
integrability in the context of Poisson manifolds.

The aim of this note is to show that the 
generalization of \cite{S,CF}
may as well be seen as a particular case of the 
previous construction in \cite{CaFeOb}. 
The main observation is that
the dual bundle $A^*$ of a Lie algebroid $A$ is naturally a Poisson manifold.
Moreover, if $A=\LIE(G)$, then $T^*G$ can be given a symplectic groupoid structure for the 
Poisson manifold $A^*$ \cite{CDW}.
As we will recall at the end of~\ref{sg}, 
it follows from classical results---\cite{M}, \cite{CW}---that
the integrability of the Lie algebroid $A$ is equivalent to the integrability of the Poisson manifold
$A^*$. We want to show that also the integration methods of \cite{CF,S} and of \cite{CaFeOb}
are equivalent. On the one hand, \cite{CaFeOb} may be recovered from \cite{CF,S} as a particular case
(by looking at the Lie algebroid of a Poisson manifold). On the other hand, we show in Sect.~\ref{sec-A*}
that \cite{CF,S} may be obtained from \cite{CaFeOb} 
by observing that to a Lie algebroid $A$ one may naturally associate a Lagrangian submanifold
of the cotangent bundle of the path space of the Poisson manifold $A^*$.
Symplectic reduction then associates to it a Lagrangian submanifold of the symplectic groupoid
$T^*G$ of $A^*$ that turns out to be (isomorphic to) the Lie groupoid $G$ of $A$.

The above construction turns out to be a special case of a more general correspondence between
coisotropic submanifolds of a Poisson manifold and Lagrangian Lie subgroupoids of its
symplectic groupoid. One direction of this correspondence (from Lagrangian subgroupoids to
coisotropic submanifolds) follows from results in \cite{X}. The other direction (the ``integration'')
is proved in Sect.~\ref{sec-coiso}.

Twisted Poisson manifolds \cite{SW} have recently become popular. 
The results of  Sect.~\ref{sec-coiso} extend to the twisted case, as we discuss
in Sect.~\ref{sec-twist}.

\begin{Ack}
I acknowledge the support of the Shapiro
fund and the
kind hospitality of Penn State University where part of this work
was completed. 
I thank Ping Xu for interesting discussions and Jim Stasheff for
useful remarks. I also thank Tudor Ratiu and Alan Weinstein
for helpful comments on infinite dimensional nightmares.
I finally thank the referee for pointing out some typos and erroneous statements
in a previous version of the paper.
\end{Ack}

\section{Preliminaries}\label{sec-prel}
In this Section we review some basic notions and fix notations.

\subsection{Lie algebroids}\label{La}
A Lie algebroid $(A,B,\rho,\Lie{\ }{\ })$ is a vector bundle $A$ over a manifold $B$
together with a bundle map $\rho\colon A\to TB$ (the {\sf anchor})
 and a Lie bracket $\Lie{\ }{\ }$ on the
real vector space $\Gamma(A)$ of sections of $A$ satisfying the following compatibility
condition:
\begin{equation}\label{e:comp}
\Lie X{fY}=f\Lie XY+\LL_{\rho_*X}f\,Y,
\qquad X,\,Y\in\Gamma(A),\ f\in C^\infty(B),
\end{equation}
where $\rho_*\colon\Gamma(A)\to\frX(B)$ is the induced map of sections and $\LL$ denotes the
Lie derivative. It follows that $\rho_*$ is a morphism of Lie algebras.
We recall some examples of Lie algebroids:
\begin{enumerate}
\item Any vector bundle with trivial anchor and Lie bracket.
\item A Lie algebra regarded as a vector bundle over a point.
\item The tangent bundle $TB$ of a manifold $B$ with the usual Lie bracket of vector fields and
$\rho$ the identity map.
\item An involutive subbundle $A$ of the tangent bundle $TB$ with the usual Lie bracket and
$\rho$ the inclusion map. 
\end{enumerate}

A Lie algebroid structure on $A$ allows one to define a differential $\delta$
on the complex $\Gamma(\Lambda^\bullet A^*)$, where $A^*$ denotes the dual bundle, by the
rules
\[
\delta f:=\rho^*\dd f,\quad f\in C^\infty(B)=\Gamma(\Lambda^0 A^*)
\]
and
\begin{multline*}
\braket{\delta\alpha}{X\wedge Y}:=-\braket\alpha{\Lie XY}
+\braket{\delta\braket\alpha X}Y
-\braket{\delta\braket\alpha Y}X,\\ 
X,\,Y\in\Gamma(A),\ \alpha\in\Gamma(A^*),
\end{multline*}
where $\rho^*\colon\Omega^1(B)\to\Gamma(A^*)$ is the transpose of $\rho_*$ and
$\braket{\ }{\ }$ is the canonical pairing of sections of $A^*$ and $A$.

If we have a bundle map $\phi$ between vector bundles $A\to B$ and $\wA\to \wB$
with base map $\varphi$,
we may define the pullback $\phi^*\colon\Gamma(\Lambda^\bullet\wA^*)\to\Gamma(\Lambda^\bullet A^*)$
as the algebra homomorphism which is $\varphi^*$ in degree zero and
the induced map\footnote{As usual, 
this is defined by setting $\phi^*(\sigma)(b)=\phi_b^*\sigma(\varphi(b))$,
$\sigma\in\Gamma(\wA^*)$, $b\in B$, where $\phi_b^*$ is the transpose of the linear map
$\phi_b\colon A_b\to\wA_{\varphi(b)}$.} 
of sections $\phi^*\colon\Gamma(\wA^*)\to\Gamma(A^*)$ in degree one.
If $A$ and $\wA$ are Lie algebroids, a bundle map $\phi$ is said to be a {\sf morphism}
if $\phi^*$ is a chain map w.r.t.\ the corresponding differentials.

If we choose local coordinates $\{b^i\}_{i=1,\dots,\dim B}$ on a trivializing
chart $U$ and
pick a basis
$\{e^\mu\}_{\mu=1,\dots,\rank A}$ on the fiber (which we also regard as a basis of
constant sections of $A_{|_U}$), 
we may introduce the anchor functions $\rho^{\mu i}$ and the structure functions $f^{\mu\nu}_\sigma$ by
the equations
\[
\rho_*(e^\mu)(b) = \rho^{\mu i}(b)\,\frac\de{\de b^i},\qquad
\Lie{e^\mu}{e^\nu}(b) = f^{\mu\nu}_\sigma(b)\, e^\sigma,
\]
where summation over repeated indices is understood. The compatibility condition \eqref{e:comp}
corresponds locally to PDEs to be satisfied by the anchor and structure functions.

\subsection{Lie groupoids}
A {\sf groupoid}
is a small category where all morphisms are invertible.
Explicitly, we have a set $G$ of morphisms and a set $B$ of objects together with
structure maps satisfying certain axioms. First we have the surjective
{\sf source} and {\sf target} maps
$s,t\colon G\to B$ and the {\sf identity} bisection $\epsilon\colon B\hookrightarrow G$. 
Then we have the multiplication $m\colon\gtwo\to G$, with
$\gtwo:=\{(u,v)\in G\times G : s(u)=t(v)\}$; as a shorthand notation we will also write $uv$ instead of
$m(u,v)$. Finally, we have an {\sf inverse} map $G\to G$, $u\mapsto u^{-1}$.
The axioms to be satisfied are 
\begin{gather*}
s(uv)=s(v),\quad t(uv)=t(u),\quad
\epsilon(b)v=v,\quad u\epsilon(b)=u,\quad
(uv)w=u(vw),\\
s(u^{-1})=t(u),\quad t(u^{-1})=s(u),\quad
uu^{-1}=\epsilon(t(u)),\quad
u^{-1}u=\epsilon(s(u)),
\end{gather*}
for all $u,v,w\in G$ and $b\in B$ for which the above expressions are meaningful.
The set $G$ is usually referred to as the groupoid, while the set $B$ is the base.
To denote a groupoid $G$ with base $B$, we will often use the notation $G\toto B$.
A {\sf morphism} between groupoids $G\toto B$ and $\wG\toto\wB$
is just a functor, i.e., a pair of maps $\Phi\colon G\to\wG$ and $\varphi\colon B\to\wB$
compatible with the structure maps $s,t,\epsilon,m$.

For a groupoid to be a {\sf Lie groupoid} \cite{P}, 
one first requires that $G$ should be a (possibly non Hausdorff) manifold, 
that $B$ should be a Hausdorff manifold and
that the source, target, identity and inverse maps should be smooth.
One further requires that the source (or equivalently the target) map should be a
submersion. This makes $\gtwo$ into a manifold, too, and one eventually requires that the
multiplication map should also be smooth. One says that a Lie groupoid is 
{\sf source-simply-connected} ({\sf ssc}) if the $s$-fibers are connected and simply connected.
A {\sf morphism} of Lie groupoids is a smooth morphism of the underlying groupoids.

Here are some examples of Lie groupoids:
\begin{enumerate}
\item $G$ a vector bundle over $B$ with $s=t=$ the projection, $\epsilon$ the zero section, multiplication
$(b,a)(b,a')=(b,a+a')$ and inverse $(b,a)^{-1}=(b,-a)$.
\item $G$ a Lie group and $B$ a point.
\item $G=B\times B$ with $s$ and $t$ the two projections, $\epsilon$ the diagonal map, multiplication
$(b,b')(b',b")=(b,b")$ and inverse $(b,b')^{-1}=(b',b)$.
\end{enumerate}

The vector bundle $\ker(\dd s)_{|_{\epsilon(B)}}$ has a natural structure of a Lie algebroid over $B$
with anchor $\dd t$ and Lie bracket induced by the multiplication. A morphism of Lie groupoids
induces a morphism of the corresponding Lie algebroids by taking its differential at the identity sections.

We will denote by
$\LIE(G)$ the Lie algebroid of the Lie groupoid $G$. The Lie algebroids of the examples
1), 2) and 3) above are the ones described in the previous subsection with the same numbers.
It is a fundamental fact that not all Lie algebroids arise in this way. Those which do are usually
called {\sf integrable}. Other fundamental facts of the theory of Lie groups generalize to Lie
groupoids:
\begin{description}
\item[Lie~I] \label{LieI} Let 
$A=\LIE(G)$ and let $\wA$ be a Lie subalgebroid of $A$. Then there is
a Lie subgroupoid of $G$ with $\wA$ as its Lie algebroid.
\item[Lie~II] Let $G$ and $\wG$ be Lie groupoids, and let $A$ and $\wA$ be the corresponding Lie algebroids.
If $G$ is ssc, for every morphism $\phi\colon A\to\wA$ there is a unique morphism 
$\Phi\colon G\to\wG$ that induces $\phi$.
\end{description}
For more information on Lie groupoids and Lie algebroids, see \cite{M}, \cite{CW} and \cite{MM}.

\subsection{Poisson manifolds}
A Poisson manifold $(M,\pi)$ is a smooth manifold $M$ together with a bivector field (i.e.,
a section of $\Lambda^2 TM$) $\pi$ such that the bracket of functions
$\Poiss fg:=\pi(\dd f,\dd g)$ satisfies the Jacobi identity.
Examples of Poisson manifolds are strong\footnote{\label{f:sympl}A symplectic form on $M$ 
is a closed, nondegenerate
$2$\ndash form. If it induces an isomorphism $TM\to T^*M$, as always happens in finite dimensions,
and not just a monomorphism,
it is called strong; otherwise it is called weak.}
symplectic manifolds (with the usual Poisson bracket of functions)
and dual spaces of Lie algebras (with linear Poisson bracket induced by the Lie bracket).
In this paper we will only consider finite-dimensional Poisson manifolds.

Lie algebroids and Poisson manifolds are strict relatives as we describe in the following.

\subsubsection{The Lie algebroid of a Poisson manifold}
To each Poisson manifold
$M$ one may associate a Lie algebroid structure on $T^*M$ with anchor the bundle map
$\pi^\#\colon T^*M\to TM$ defined by  $\pi^\#(x)(\sigma)=\pi(x)(\sigma,\bullet)$, with $x\in M$,
$\sigma\in T^*_xM$, and the {\sf Koszul bracket} which on
exact forms is defined by $\Lie{\dd f}{\dd g}:=\dd\Poiss fg$ and is extended to arbitrary
$1$\ndash forms using the rule \eqref{e:comp}.\footnote{\label{foot:delta}The corresponding 
differential on sections of the exterior algebra of the dual of $T^*_\pi M$
(i.e., on multivector
fields) is the inner derivation $\delta:=\Lie\pi{\ }_\text{SN}$ where $\Lie{\ }{\ }_\text{SN}$ 
is the Schouten--Nijenhuis bracket.
So $(T^*_\pi M,TM)$ is an example of Lie bialgebroid.}
We will denote the Lie algebroid of the Poisson manifold
$(M,\pi)$ by $T^*_\pi M$. One says that a Poisson manifold is {\sf integrable} if its
Lie algebroid is.

In local coordinates $\{x^i\}_{i=1,\dots,\dim M}$,
we have the local basis ${\dd x^i}$ of constant $1$\ndash forms w.r.t.\ which
the anchor functions are the components $\pi^{ij}$ and the structure functions are the partial
derivatives $\de_i\pi^{rs}$.

\subsubsection{Lie algebroids as Poisson manifolds}\label{A*}
If $(A,B,\rho,\Lie{\ }{\ })$ is a Lie algebroid,
the dual bundle $A^*$ has a natural Poisson structure.
This is defined first on functions that are constant on the fibers 
(i.e., functions on $B$) or linear on the fibers (i.e., sections
of $A$):
\[
\Poiss FG = 
\begin{cases}
0 &\text{if $F,G\in C^\infty(B)$,}\\
\LL_{\rho_*F}G &\text{if $F\in\Gamma(A)$ and $G\in C^\infty(B)$,}\\
\Lie FG &\text{if $F,G\in\Gamma(A)$}.
\end{cases}
\]
The bracket is then
extended to functions that are polynomial on the fibers (i.e., sections of the symmetric
powers of $A$) 
as a skew-symmetric biderivation
and finally to all smooth functions by completion.

If we choose local coordinates $\{b^i\}_{i=1,\dots,\dim B}$ and $\{\alpha^\mu\}_{\mu=1,\dots,\rank A}$ on
$A^*$, then the local coordinate expression of the bivector field corresponding to the above Poisson 
structure is
\begin{equation}\label{piA}
\pi(b,\alpha) =  \alpha^\sigma f^{\mu\nu}_\sigma(b)\frac\de{\de\alpha^\mu}\wedge\frac\de{\de\alpha^\nu} +
\rho^{\mu i}(b)\frac\de{\de\alpha^\mu}\wedge\frac\de{\de b^i},
\end{equation}
where $\rho^{\mu i}$ and $f^{\mu\nu}_\sigma$ are the anchor and
structure functions respectively.

\subsection{Symplectic groupoids}\label{sg}
If $G$ is a Lie groupoid,
there are three natural maps from $\gtwo$ to $G$: the two projections $p_1$ and $p_2$ and the
multiplication $m$. A differential form $\omega$ on $G$ is said to be {\sf multiplicative} if it
satisfies the cocycle condition
\[
m^*\omega = p_1^*\omega+p_2^*\omega.
\]
(If $\omega$ is a function, this simply means $\omega(uv)=\omega(u)+\omega(v)$ for all $u,v\in\gtwo$.)

A {\sf symplectic groupoid} \cite{CDW,K,Z}
$(G\toto M,\omega)$
is then by definition a Lie groupoid $G$ over $M$ endowed with a multiplicative symplectic
form $\omega$. It follows \cite{CDW} that $\epsilon\colon M\hookrightarrow G$ is a Lagrangian
embedding, that the inverse map is an anti-symplectomorphism, and that
the base manifold $M$ has a unique Poisson structure
$\pi$ such that the source and the target maps are Poisson and anti-Poisson respectively.
The Lie algebroid of $(G\toto M,\omega)$ turns then out to be isomorphic to $T^*_\pi M$.

It is proved in \cite{MX} that, given an integrable Poisson manifold $(M,\pi)$, it is alway possible
to endow a ssc Lie groupoid $G$ such that $\LIE(G)=T^*_\pi M$ with a multiplicative symplectic
form such that the induced Poisson structure on $M$ is $\pi$.
This means that the problem of integrating a Poisson manifold to its symplectic groupoid
actually amounts just to integrating its Lie algebroid.

We discuss one single example of symplectic groupoid, which is relevant for the rest of the paper.
Let $G\toto B$ be a Lie groupoid. Endow $T^*G$ with the canonical symplectic structure $\omega$.
Let $A=\LIE(G)$; then one can endow $T^*G$ with a Lie groupoid structure with base $A^*$
such that $\omega$ is multiplicative and the Poisson structure on $A^*$ is the one described 
in~\ref{A*} (see \cite{CDW}). Moreover, if $G$ is ssc, then so is $T^*G$. 
This shows that $A^*$ is integrable as a Poisson manifold if $A$ is integrable as a Lie algebroid.
The converse is also true by Lie~I on page~\pageref{LieI} since $A$ may be regarded as a Lie subalgebroid
of $T^*A^*$ (see Lemma~\ref{lem:AA}).

Thus, in order to integrate a Lie algebroid $A$, one may
first look for the ssc symplectic groupoid of $A^*$ and the look for the Lagrangian Lie subgroupoid
whose Lie algebroid is $A$. This is not immediate as the ssc symplectic groupoid of $A^*$ 
may not be presented as $T^*G$. In Sect.~\ref{sec-A*} we will however describe
a method to do this explicitly.

\subsection{Symplectic reduction}\label{sr}
Let $(S,\omega)$ be a (possibly weak) symplectic manifold (see footnote~\ref{f:sympl})
and $C$ a submanifold.
The orthogonal tangent bundle $T^\perp C$ is defined as the subbundle of $T_CS$ consisting
of vectors that are $\omega$\ndash orthogonal to vectors tangent to $C$.
The submanifold $C$ is called {\sf coisotropic} if $T^\perp C\subset TC$ and---as a particular 
case---{\sf Lagrangian} if $T^\perp C=TC$. Since $\omega$ is closed, the subbundle $T^\perp C$
defines an involutive distribution (the {\sf characteristic foliation})
on the coisotropic submanifold $C$ whose leaves
are exactly the kernel of $\omega$. This implies that the leaf space $\underline C$, if smooth,
is naturally endowed with a symplectic form $\underline\omega$ whose pullback to $C$ is the restriction
of $\omega$.

If $C$ is coisotropic, $L$ is Lagrangian and their intersection is {\sf clean}
(viz., $L\cap C$ is also a submanifold and $T(L\cap C)=TL\cap TC$),
then the image of the projection of $L\cap C$ to $\underline C$ is Lagrangian in $\underline C$
(it may not be an embedding but just an immersion), see \cite[Lecture~3]{W}.
For this to hold in infinite dimensions, one has to add explicitly a further condition
(which in finite dimensions is automatically satisfied):
$T^\perp(L\cap C)=TL+T^\perp C$. When this condition is satisfied we say that the intersection
is {\sf symplectically regular}.

Even when $\underline C$ is not smooth, we may think of it as a singular symplectic manifold
and of the projection of $L\cap C$ as a singular Lagrangian submanifold.

\section{Integration of Poisson manifolds}\label{sec-ipm}
We briefly recall the method of \cite{CaFeOb} and fix the notations.
Let $M$ be a finite-dimensional
Poisson manifold with Poisson bivector $\pi$.
Let $PM:=\{I\to M\}$ be the path space of $M$ and $T^*PM$ the manifold\footnote{A Banach manifold
structure may be introduced by restricting to maps with a given degree of differentiability.
For example, as in \cite{CaFeOb}, one may define $T^*PM$ as the space of continuous bundle maps with
$C^1$\ndash base maps.}
of bundle
maps $TI\to T^*M$. The fiber over $X\in PM$ may be identified with the space of sections 
$\Gamma(T^*I\otimes X^*T^*M)$.
The canonical symplectic form $\Omega$ at a point $X\in PM$, $\eta\in T^*_X PM$ is defined by
\begin{multline*}
\Omega{(X,\eta)}(\xi_1\oplus e_1,\xi_2\oplus e_2)=
\int_I \braket{e_1}{\xi_2} - \braket{e_2}{\xi_1},\\
\xi_i\in \Gamma(X^*TM)=T_XPM,\quad
e_i\in \Gamma(T^*I\otimes X^*T^*M)=T_X^*PM,
\end{multline*}
where $\braket{\ }{\ }$ denotes the canonical pairing between
tangent and cotangent fibers of $M$.
The submanifold $\calC(M)$ defined as the space of solutions to the 
equations\footnote{As observed in \cite{S,CF}, solutions to \eqref{C}---called
``cotangent paths'' in \cite{CFPoiss}---are precisely Lie algebroid morphisms
$TI\to T^*_\pi M$, where $TI$ is given its canonical Lie algebroid structure.}
\begin{equation}\label{C}
\dd X=\pi^\#(X)\eta
\end{equation}
is coisotropic 
in $T^*PM$.
Its characteristic foliation turns out to be expressed in terms of an infinitesimal action of the Lie algebra
\[
\Gamma_0(M) := \{C\colon I\to \Omega^1(M)\ |\ C(0)=C(1)=0\}, 
\]
where the Lie bracket is defined pointwise in terms of the Lie bracket on $\Omega^1(M)$.
To describe this action, it is easier to pass to local coordinates
$\{x^i\}_{i=1,\dots,\dim M}$ on $M$. Then the vector field on $\calC(X)$ associated to an element $C$
of $\Gamma_0(M)$ evaluated at a point $(X,\eta)$ may be written as $(\delta X^i,\delta \eta_i)$
with
\begin{subequations}\label{delta}
\begin{align}
\delta X^i &= -\pi^{ij}(X)\,(C_X)_j,\\
\delta\eta_i &= -\dd (C_X)_i - \de_i\pi^{rs}\,(C_X)_r\,\eta_s,
\end{align}
\end{subequations}
where $C_X$ is the section of $X^*T^*M$ defined by $C_X(t)=C(t)(X(t))$ and $\de_i\pi^{rs}$
are the structure functions of the Lie algebroid $T^*_\pi M$ 
(w.r.t.\ the local basis $\{\dd x^i\}$ of sections).
Upon using the constraint equations \eqref{C},
this local coordinate expression is well-defined. The leaf space 
$\underline\calC(M)$
of $\calC(M)$ is then
the (possibly singular) ssc symplectic groupoid of $M$.

\section{Integration of $A$ and $A^*$}\label{sec-A*}
We now apply the method recalled in Sect.~\ref{sec-ipm} to the Poisson manifold $A^*$, where
$A$ is a Lie algebroid.
We denote by $\calP(A)$ the manifold of bundle maps $TI\to A$.
The central result of this paper is the following 
\begin{Thm}\label{thm-A*}
$\calP(A)$ is a Lagrangian submanifold of $T^*PA^*$, and the projection $\G(A)$
of $\calP(A)\cap\calC(A^*)$ to
$\underline\calC(A^*)$
yields the (possibly singular) ssc Lie groupoid of $A$ as a Lagrangian Lie subgroupoid
of the symplectic groupoid of $A^*$.
\end{Thm}
The rest of the Section is devoted to the proof of this Theorem.

\begin{proof}
We begin with an easy Lemma whose proof is left to the reader.
\begin{Lem}\label{lem:AA}
Let $A\to B$ be a vector bundle. 
The fiber of $T^*A^*$ over a point $(b,\alpha)\in A$
is the vector space $T^*_b B\oplus A_b$. 
The map 
\[
\iota\colon
\begin{array}[t]{ccc}
A &\to & T^*A^*\\
(b,a) &\mapsto & ((b,0),0\oplus a)
\end{array}
\]
is an injective bundle map from $\begin{matrix}A\\ \downarrow\\ B\end{matrix}$ to 
$\begin{matrix}T^*A^*\\ \downarrow\\ A^*\end{matrix}$. 

If $T^*A^*$ is given
the canonical symplectic structure, then $\iota$ is a Lagrangian embedding.

If $A$ is a Lie algebroid and $T^*A^*$ is given the Lie algebroid structure induced by the Poisson
structure on $A^*$, then $\iota$ is a morphism of Lie algebroids. So $A$ is a Lagrangian Lie subalgebroid
of $T^*A^*$.
\end{Lem}
As a consequence, the composition of a bundle map $TI\to A$ with $\iota$ yields a bundle map
$TI\to T^*A^*$. So $\iota$ induces an inclusion of $\calP(A)$ into $T^*PA^*$ which is also Lagrangian.

\begin{Lem}
The intersection $\calP(A)\cap\calC(A^*)$ consists of Lie algebroid morphisms $TI\to A$.
\end{Lem}
\begin{proof}
Choosing local coordinates as in Sect.~\ref{La}, we denote an element of $T^*PA^*$ by the
functions $X^i$, $\alpha^\mu$ together with the $1$\ndash forms $\eta_i$ and $a_\mu$
(observe that at different
points in $I$ we may be on different patches of local coordinates).
Using the Poisson bivector field $\pi$ defined in \eqref{piA},
the constraint equations \eqref{C} read
\begin{subequations}
\begin{align}
\dd X^i &= \rho^{\mu i}(X) a_\mu,\label{se:dXa}\\
\dd\alpha^\mu &= -\rho^{\mu i}(X) \eta_i - \alpha^\sigma f^{\mu\nu}_\sigma(X) a_\nu.\label{se:daeta}
\end{align}
\end{subequations}
Elements of $\calP(A)$ are represented by functions and $1$\ndash forms as above with the conditions
$\alpha^\mu=0$ and $\eta_i=0$ $\forall \mu,i$. So the intersection consists
of functions  $X^i$ together with $1$\ndash forms $a_\mu$ satisfying
the first equation above, which is precisely the local coordinate expression of the condition
that the bundle map $TI\to A$ defined by $(X,a)$ is a morphism of Lie algebroids.
Observe finally that both sides of \eqref{se:daeta} identically vanish, so no further conditions are imposed.
\end{proof}
In the following we will denote by $\PP(A)$
the space of  Lie algebroid morphisms $TI\to A$.
Given a second unit interval $J=[0,1]$, we denote by $\PP_2(A)$
the space of  Lie algebroid morphisms $T(I\times J)\to A$.
Two elements $\gamma_0$ and $\gamma_1$ of $\PP(A)$ are said to be
Lie algebroid homotopic, if there exists an element
of $\PP_2(A)$ that restricts to $\gamma_u$ at $TI\times\{u\}$, $u=0,1$,
and is trivial at $\{v\}\times TJ$, $v\in\de I$.

\begin{Lem}
The foliation of $\calC(A^*)$ restricted to $\calP(A)\cap\calC(A^*)$ is precisely
the infinitesimal version of Lie algebroid homotopies.
\end{Lem}
\begin{proof}
The foliation is defined by \eqref{delta}.
Given a map $X\colon I\to B$,
an element $C$ of the Lie algebra $\Gamma_0(A^*)$ yields a section $C_X$ of $X^*T^*A^*\to I$.
Again choosing local coordinates, we may denote the components of $C_X(t)$ in $T^*_{X(t)}B$ by
$\beta_i(t)$ and those in $A_{X(t)}$ by $b_\mu(t)$, $t\in I$.
The vector field corresponding to $(\beta,b)$ 
evaluated at a point $(X,\alpha,\eta,a)$
may be written using
\eqref{delta} and \eqref{piA} as $(\delta X^i,\delta\alpha^\mu,\delta\eta_i,\delta a_\mu)$ with 
\begin{align*}
\delta X^i &= -\rho^{\mu i}(X) b_\mu,\\
\delta\alpha^\mu &= \rho^{\mu i}(X)\beta_i + \alpha^\sigma f^{\mu\nu}_\sigma(X) b_\nu,\\
\delta\eta_i &= -\dd\beta_i - \alpha^\sigma \de_if^{\mu\nu}_\sigma(X)a_\mu b_\nu 
-\de_i\rho^{\mu j}(X)(a_\mu\beta_j-\eta_jb_\mu),\\
\delta a_\mu &= -\dd b_\mu - f^{\nu\sigma}_\mu(X)a_\nu b_\sigma.
\end{align*}
The restriction of this foliation to $\calP(A)\cap\calC(A^*)$  is given by the first and last equations
(the remaining ones are automatically satisfied when we impose $\alpha=\eta=0$ and, consequently,
$\delta\alpha=\delta\eta=0$).
If we integrate the flow of this vector field on a time-interval $(-\epsilon,\epsilon)$, 
these equations are then 
precisely the local coordinate expressions of a morphism of Lie algebroids
$T(I\times(-\epsilon,\epsilon))\to A$.
\end{proof}
This shows that the projection of $\calP(A)\cap\calC(A^*)$ to
$\underline\calC(A^*)$ is equal to the quotient $\G(A)$ 
of $\PP(A)$ by Lie algebroid morphisms.
It was shown in \cite{CF}, along the lines of
\cite{CaFeOb}, that $\G(A)$
has a groupoid structure and that, if it is smooth, has
$A$ as its Lie algebroid.

\begin{Lem}
The intersection of $\calP(A)$ and $\calC(A^*)$ is clean and symplectically regular.
\end{Lem}
\begin{proof}
The intersection $\calP(A)\cap\calC(A^*)$ may be given a manifold structure, see \cite{CF}.
An element $(X,a)$ in it is a solution to \eqref{se:dXa}.
So a tangent vector at $(X,a)$ is a pair $(\dot X,\dot a)$ satisfying
\begin{equation}\label{e:ddotX}
\dd \dot X^i = \rho^{\mu i}(X) \dot a_\mu
+ \dot X^j\de_j\rho^{\mu i}(X) a_\mu.
\end{equation}
On the other hand, the tangent vector at $(X,\alpha,\eta,a)\in\calC(A^*)$
consists of a quadruple $(\dot X,\dot \alpha,\dot \eta,\dot a)$
satisfying
\begin{align*}
\dd \dot X^i &= \rho^{\mu i}(X) \dot a_\mu +
\dot X^j\de_j\rho^{\mu i}(X) a_\mu,\\
\dd\dot\alpha^\mu &= -\rho^{\mu i}(X) \dot\eta_i 
-\dot X^j\de_j\rho^{\mu i}(X) \eta_i +\\
 &\phantom{=}
\ - \alpha^\sigma f^{\mu\nu}_\sigma(X) \dot a_\nu
- \dot \alpha^\sigma f^{\mu\nu}_\sigma(X) a_\nu
- \alpha^\sigma \dot X^j\de_jf^{\mu\nu}_\sigma(X) a_\nu.
\end{align*}
At an intersection point with $\calP(A)$
we have to set $\alpha=0$ and $\eta=0$, and intersecting with $T\calP(A)$ means setting
$\dot\alpha=0$ and $\dot\eta=0$; so we recover \eqref{e:ddotX}. This shows that the intersection
in clean.

Elements of $T\calP(A)+T\calC(A^*)$ at a point $(X,0,0,a)\in\calP(A)\cap\calC(A^*)$
are vectors of the form $(\dot X,\delta\alpha,\delta\eta,\dot a)$ where
$\dot X$ and $\dot a$ are arbitrary while
\begin{align*}
\delta\alpha^\mu &= \rho^{\mu i}(X)\beta_i,\\
\delta\eta_i &= -\dd\beta_i  
-\de_i\rho^{\mu j}(X)a_\mu\beta_j.
\end{align*}
An explicit, though lengthy, computation shows that these are precisely all possible vectors
in $T^\perp(\calP(A)\cap\calC(A^*))$.
\end{proof}
By the discussion in~\ref{sr} 
we conclude the proof of the Theorem.
\end{proof}

If we recall that the Lie groupoid of $A$ must appear as a Lagrangian Lie subgroupoid of the
symplectic groupoid of $A^*$ (see the end of~\ref{sg}), 
we may interpret the above result as a way of deriving
the method of \cite{S,CF} from the one of \cite{CaFeOb}.

Finally observe that, with the above notation, $\calC(M)=\PP(T^*_\pi M)$
and $\underline\calC(M)=\G(T^*_\pi M)$, 
so the method
of \cite{CaFeOb} may also be recovered as a particular case of \cite{S,CF}.

\section{Integration of coisotropic submanifolds}\label{sec-coiso}
Let $(M,\pi)$ be a Poisson manifold. A submanifold $C$ is called {\sf coisotropic}
if $\pi^\#(N^*C)\subset TC$, where $N^*C$ denotes the conormal bundle of $C$
(viz., the subbundle of $T^*_CM$ of covectors that vanish when applied to a vector tangent to $C$).
In case the Poisson structure of $M$ comes from a strong symplectic structure, this definition
coincides with the usual one recalled in~\ref{sr} since, in this case, $\pi^\#$ establishes an isomorphism
between $N^*C$ and $T^\perp C$.

The theory of coisotropic submanifolds of Poisson manifolds \cite{W2}
generalizes many properties of the corresponding
theory in the symplectic case (e.g., one may generalize symplectic reduction to the Poisson case
as it turns out that $\pi^\#(N^*C)$ is an integrable distribution on the coisotropic submanifold $C$
and that the leaf space inherits a Poisson structure).
Moreover, coisotropic submanifolds label the possible boundary conditions of the Poisson sigma model
yielding the beginning of a theory of quantum reduction in the deformation quantization context
\cite{CFcoiso}.

{}From the point of view of the present paper, coisotropic submanifolds are important because of their
relations with the theory of Lie algebroids. Namely:
\begin{Prop}\label{coisosub}
A submanifold $C$ of $(M,\pi)$ is coisotropic if{f} $N^*C$ is a Lie subalgebroid of $T^*_\pi M$.
\end{Prop}
\begin{Prop}\label{coisoA}
Let $A$ be a Lie algebroid and $A^*$ its dual regarded as a Poisson manifold.
Then the zero section of $A^*$ is coisotropic and its conormal bundle
is the inclusion $\iota$ of $A$ as a Lagrangian Lie subalgebroid
of $T^*A^*$ described in Lemma~\ref{lem:AA}.
\end{Prop}
\begin{proof}[Proof of Prop.~\ref{coisosub}]
If $N^*C$ is a Lie subalgebroid, in particular its anchor $N^*C\to TC$ is the restriction of $\pi^\#$
to $N^*C$; this immediately shows that $C$ is coisotropic.

The converse is true by Corollary 3.1.5 in \cite{W2}, but for completeness we give a proof here.
Assume that $C$ is coisotropic. By definition, the restriction of $\pi^\#$
to $N^*C$ maps it to $TC$ and so it defines an anchor. It remains only to prove that
the Koszul bracket induces a bracket on sections of $N^*C$.
Let $U$ be a trivializing chart on $M$ intersecting
$C$. We choose adapted local coordinates $\{x^I\}_{I=1,\dots,\dim M}$
so that $U\cap C$ is determined by $x^I=0$, $I=\dim C+1,\dots,\dim M$. 
To make the notation more transparent, we will use small Latin indices to denote the first
$\dim C$ coordinates (the tangential ones) and small Greek indices to denote the remaining 
(transversal) ones; when we do not want to distinguish them, we will use capital Latin indices.
So the above conditions may be written $x^\mu=0$.
If $C$ is coisotropic, then $\pi^{\mu\nu}(x)=0$ for all $x\in U\cap C$, and as a consequence
$\de_i\pi^{\mu\nu}(x)=0$ for all $x\in U\cap C$. This implies that
\[
\Lie{\dd x^\mu}{\dd x^\nu}(x)=
\de_K\pi^{\mu\nu}(x)\,\dd x^K =
\de_\mu\pi^{\mu\nu}(x)\,\dd x^\mu,
\qquad\forall x\in U\cap C.
\]
This concludes the proof.
\end{proof}
\begin{proof}[Proof of Prop.~\ref{coisoA}]
Using coordinates $(b,\alpha)$
as in~\ref{A*}, we see that on the zero section $\alpha=0$ only the mixed components
of the bivector field in \eqref{piA} survive. This shows that the zero section is coisotropic.
Its conormal bundle consists of elements in $T^*A^*$ of the form $((b,0),0\oplus a)$, so
it is $\iota(A)$.
\end{proof}

Observe that the conormal bundles of submanifolds of $M$ are all possible Lagrangian subbundles
of $T^*M$ with its canonical symplectic structure. So Prop.~\ref{coisosub} may also be rephrased as
\begin{Prop}\label{coLa}
The set of coisotropic submanifolds of $M$ is isomorphic to the set of Lagrangian Lie subalgebroids
of $T^*_\pi M$. 
\end{Prop}

In \cite[Sect.~4]{X} coisotropic subgroupoids of Poisson--Lie groupoids are studied.
It follows from Prop.~4.10 there, as a particular case, that the Lie algebroid of a Lagrangian Lie subgroupoid
of the symplectic groupoid of the Poisson manifold $(M,\pi)$ is a Lagrangian Lie subalgebroid of
$T^*_\pi M$. We want to prove that also the converse of this statement is true.
\begin{Thm}\label{Dumbledore}
Let $(M,\pi)$ be an integrable 
Poisson manifold and $\underline\calC(M)$ its ssc symplectic groupoid.
Then there is a one-to-one correspondence between Lagrangian Lie subgroupoids of 
$\underline\calC(M)$ and coisotropic submanifolds of $M$.
\end{Thm}
One direction of the isomorphism follows from the cited Prop.~4.10 of \cite{X} 
together with Prop.~\ref{coLa} above.
We have then to construct an inverse map from coisotropic submanifolds of $M$ to 
Lagrangian Lie subgroupoids of  $\underline\calC(M)$. We do it using the technique
of \cite{CaFeOb} recalled in
Sect.~\ref{sec-ipm} and actually prove a more general statement:
\begin{Prop}\label{Potter}
To each coisotropic submanifold $C$ of $M$ there corresponds a (possibly singular)
Lagrangian Lie subgroupoid (isomorphic to $\G(N^*C)$ as a groupoid) of $\underline\calC(M)$.
\end{Prop}
Observe then that, thanks to Prop.~\ref{coisoA},
Thm.~\ref{thm-A*} is now a particular case of this Proposition.
To prove the Proposition we introduce 
\begin{equation}\label{e:calLC}
\calL(C):=\{(X,\eta)\in T^*PM : X\in PC,\ \eta\in\Gamma(T^*I\otimes X^*N^*C)\}.
\end{equation}
It is easy to see that $\calL(C)$ is a Lagrangian submanifold of $T^*PM$ and
that $\calL(C)\cap\calC(C)$ is the manifold $\PP(N^*C)$ of Lie algebroid morphisms
$TI\to N^*C$. To complete the proof of Prop.~\ref{Potter} (and hence of Thm.~\ref{Dumbledore}),
by the discussion in~\ref{sr}
we only need the following two Lemmata:
\begin{Lem}\label{LC}
The restriction to $\calL(C)\cap\calC(C)$ of the characteristic foliation of $\calC(C)$
is the infinitesimal form of Lie algebroid homotopies $\PP_2(N^*C)$.
\end{Lem}
\begin{Lem}\label{coisoclean}
The intersection $\calL(C)\cap\calC(C)$ is clean and symplectically regular.
\end{Lem}
\begin{proof}[Proof of Lemma~\ref{LC}]
For simplicity we use
adapted local coordinates as in the proof of Prop.~\ref{coisosub} (again observe that at different
points in $I$ we may be on different patches of adapted local coordinates).
The intersection of the foliation \eqref{delta} with $\calL(C)$ amounts to
the constraints $\delta X^\mu=0$ and $\delta\eta_i=0$.
On $\calL(C)$ we have $\eta_i=0$ and $\pi^{\mu\nu}(X)=0$; so
\[
\delta\eta_i = -\dd (C_X)_i - \de_i\pi^{RS}\,\eta_R\,(C_X)_S=
-\dd (C_X)_i - \de_i\pi^{\mu k}\,\eta_\mu\,(C_X)_k,
\]
and
the condition $\delta\eta_i=0$ together with the boundary
conditions on $C_X$ implies $(C_X)_i=0$. Then $\delta X^\mu = -\pi^{\mu j}(X)\,(C_X)_j$
automatically vanishes. Moreover, we get
\begin{align}
\delta X^i &= -\pi^{i\mu}(X)\,(C_X)_\mu,\\
\delta\eta_\mu &= -\dd (C_X)_\mu - \de_\mu\pi^{\nu\tau}\,(C_X)_\nu\,\eta_\tau.
\end{align}
The local flow of this vector field on a time interval $(-\epsilon,\epsilon)$ precisely
defines a Lie algebroid morphism $T(I\times(-\epsilon,\epsilon))\to N^*C$.
\end{proof}
\begin{proof}[Proof of Lemma~\ref{coisoclean}]
As already recalled $\calL(C)\cap\calC(C)=\PP(N^*C)$ is a manifold, see \cite{CF}.
Let now $(X,\eta)\in\calL(C)$. Using again adapted local coordinates, we see that
a vector $(\dot X,\dot\eta)\in T_{(X,\eta)}T^*PM$ belongs to 
$T_{(X,\eta)}\calC(M)$ if{f}
\[
\dd\dot X^I = \pi^{IK}(X)\,\dot\eta_K + \dot X^J\de_J\pi^{I\nu}(X)\,\eta_\nu,
\]
where we have used $\eta_i=0$. The intersection of $T_{(X,\eta)}\calL(C)$ with
$T_{(X,\eta)}\calC(M)$ is then determined by also imposing the equations 
$\dot X^\mu=\dot\eta_i=0$ (which express belonging to $T_{(X,\eta)}\calL(C)$); viz.,
\begin{align*}
\dd\dot X^i &= \pi^{i\nu}(X)\,\dot\eta_\nu + \dot X^j\de_j\pi^{i\nu}(X)\,\eta_\nu,\\
\dd\dot X^\mu &= \pi^{\mu\nu}(X)\,\dot\eta_\nu + \dot X^j\de_j\pi^{\mu\nu}(X)\,\eta_\nu.
\end{align*}
The first equation says that $(\dot X^i,\dot\eta_\mu)$ belongs to $T_{(X,\eta)}(\calL(C)\cap\calC(C))$.
The condition $\dot X^\mu=0$ implies $\dd \dot X^\mu=0$, but thanks to
$\pi^{\mu\nu}(X)=0$ and  $\de_j\pi^{\mu\nu}(X)=0$ the second equation set to zero does not put
extra conditions on $(\dot X^i,\dot\eta_\mu)$. This shows that the intersection is clean.

An explicit but lengthy computation shows that the intersection is also symplectically
regular and in particular that the elements of
\[
T_{(X,\eta)}(\calL(C)\cap\calC(C))=T_{(X,\eta)}\calL(C)+T_{(X,\eta)}\calC(C),
\]
$(X,\eta)\in\calL(C)\cap\calC(C)$, are of the form $(\dot X^i,\delta X^\mu,\delta\eta_i,\dot\eta_\mu)$
with
\begin{align*}
\delta X^\mu &= -\pi^{\mu j}(X)\,(C_X)_j,\\
\delta\eta_i &= -\dd (C_X)_i - \de_i\pi^{r\nu}\,(C_X)_r\,\eta_\nu,
\end{align*}
and $\dot X^i$ and $\dot\eta_\mu$ arbitrary.
\end{proof}

\section{The twisted case}\label{sec-twist}
A {\sf twisted symplectic manifold} is a manifold endowed with
a nondegenerate $2$\ndash form. A {\sf twisted Poisson manifold} $(M,\pi,\phi)$ is a manifold
$M$ endowed with a bivector field $\pi$ and a closed $3$\ndash form $\phi$ such that
\[
\Lie\pi\pi=\frac12 \wedge^3\pi^\#\phi.
\]
One may still define a bracket on functions by $\Poiss fg=\pi(\dd f, \dd g)$, but the Jacobi
identity will not be satisfied.
A twisted strong symplectic manifold $(M,\omega)$ provides an example of twisted Poisson manifold by
setting $\pi$ to be the inverse of $\omega$ and $\phi=\dd\omega$.
The cotangent bundle of a twisted Poisson manifold is a Lie algebroid with $\pi^\#$ as its anchor
and a Lie bracket that on exact $1$\ndash forms reads 
$\Lie{\dd f}{\dd g}=\dd\Poiss fg +\iota_{\pi^\#\dd f}\iota_{\pi^\#\dd g}\phi$.
We will denote this Lie algebroid by $T^*_{\pi,\phi}M$.

A {\sf twisted symplectic groupoid} $(G\toto M,\omega,\phi)$ is a Lie groupoid $G\toto M$ endowed
with a nondegenerate, multiplicative $2$\ndash form $\omega$ on $G$ and with a closed $3$\ndash form $\phi$
on $M$ such that the cocycle condition $\dd\omega=s^*\phi-t^*\phi$ holds.
It turns out \cite{CX} that a twisted Poisson structure $(\pi,\phi)$ is induced on the base $M$
such that $\LIE(G)$ is isomorphic to $T^*_{\pi,\phi}M$.

{\sf Symplectic reduction} may be generalized to the twisted case. Let $(S,\omega)$ be a twisted symplectic
manifold. A submanifold $C$ is called {\sf coisotropic} if $T^\perp C\subset TC$ and
the restriction $\omega_C$
of $\omega$ to $C$ is invariant (viz., $\LL_X \omega_C=0$ for any 
$X\in\Gamma(T^\perp C)$). A {\sf Lagrangian} submanifold is a submanifold $L$ such that
$T^\perp L= TL$, and it is automatically coisotropic.
It turns out that $T^\perp C$ defines a foliation on the coisotropic submanifold $C$
and that the leaf space $\underline C$ inherits a twisted symplectic structure (if smooth).
Moreover, if $L$ is Lagrangian and the intersection $L\cap C$ is clean and symplectically regular, 
the image of the projection of $L\cap C$
to $\underline C$ is also Lagrangian.

We now recall the method introduced in \cite{CX} to integrate twisted Poisson manifolds to (possibly
singular) twisted symplectic groupoids by modifying the method of \cite{CaFeOb}.
One considers again the submanifold $\calC(M)$ of $T^*PM$ as in Sect.~\ref{sec-ipm}. This is not coisotropic
w.r.t.\ the canonical symplectic form $\Omega$ of $T^*PM$, but it is so w.r.t.\ the twisted
symplectic form $\wOmega:=\Omega+\hOmega$ with
\[
\hOmega(X,\eta)(\xi_1\oplus e_1,\xi_2\oplus e_2)=
\frac12 \int_I \phi(X)(\pish(X)\eta,\xi_1,\xi_2),
\]
$\xi_1,\xi_2\in\Gamma(X^*TM)$, $e_1,e_2\in \Gamma(T^*I\otimes X^*T^*M)$.
The twisted symplectic groupoid of $M$ turns then out to be the leaf space $\underline \calC(M)$.

Finally, we want to generalize the results of Sect.~\ref{sec-coiso}.
We say that a submanifold $C$ is {\sf coisotropic} in the twisted Poisson manifold  $(M,\pi,\phi)$
if $\pi^\#(N^*C)\subset TC$ and the restriction $\phi_C$
of $\phi$ to $C$ is horizontal (viz.,
$\iota_X\phi_C=0$ for any $X\in\Gamma(N^*C)$).
If the twisted Poisson structure comes from a twisted strong symplectic structure, then this definition
coincides with the one given in the twisted symplectic case.

One may easily see (along the lines of Sect.~\ref{sec-coiso}) that, if $C$ is coisotropic,
$N^*C$ is a Lie subalgebroid of $T^*_{\pi,\phi}M$ and that 
the leaf space $\underline C$ inherits the structure of a twisted Poisson manifold.
However, $N^*C$ might be a Lie subalgebroid of $T^*_{\pi,\phi}M$ in more general instances.
Thus, we no longer have a one-to-one correspondence between coisotropic submanifolds
of $(M,\pi,\phi)$ and Lagrangian Lie subalgebroids of $T^*_{\pi,\phi}M$ if on $T^*M$ we put 
the canonical symplectic structure
$\omega=\dd p_i\,\dd x^i$. 
On the other hand, if we put on $T^*M$ the
twisted symplectic structure\footnote{In \cite{CX}
it is observed that $(T^*_{\pi,\phi}M,TM)$ is actually an example of quasi Lie bialgebroid
(cf.\ footnote~\fullref{foot:delta})
as there is a derivation $\delta$ of $\Omega^\bullet(M)$ (which deforms the exterior derivative)
such that 
\[
\delta\Lie\sigma\tau=\Lie{\delta\sigma}\tau+\Lie\sigma{\delta\tau},
\qquad\forall\sigma,\tau\in\Omega^1(M),
\]
and that $\delta^2=\Lie\phi{\bullet}$
(where we have extended the Lie bracket to the whole of $\Omega^\bullet(M)$
as a biderivation).
This corresponds to having a quasi (i.e., no Jacobi)
Lie algebroid structure on $TM$
that is compatible with the one on $T^*_{\pi,\phi}M$. 
This quasi Lie algebroid structure determines a twisted Poisson structure on $T^*M$, and as this
turns out to be nondegenerate, it corresponds to a twisted symplectic structure that is precisely
$\womega$.} 
$\womega:=\omega+\homega$
with $2\,\homega=p_i\pi^{ij}\phi_{jkl}\,\dd x^k\,\dd x^l$,
it turns out that $N^*C$ is Lagrangian in $(T^*M,\womega)$
if{f} $\phi_C$ is horizontal; so
Prop.~\ref{coLa} generalizes to the twisted case, provided one twists the symplectic form
on $T^*M$.
One may finally generalize Prop.~\ref{Potter} to the twisted case since
$\calL(C)$, defined as in \eqref{e:calLC}, turns out to be Lagrangian
in $(T^*PM,\wOmega)$ if $C$ is coisotropic, and
it is possible to show Lemmata~\ref{LC} and~\ref{coisoclean} in this case
(with very similar proofs).

A further modification of the integration methods of \cite{CaFeOb, CX} is considered in \cite{BCWZ}
to integrate (twisted) Dirac structures, a common
generalization of (twisted) symplectic and Poisson structures.
It would be interesting to know if any of the ideas in the present paper have a generalization
in that context.

\end{document}